\documentstyle{article}
\begin{document}
\title{{On Axially Symmetric  Solutions of Fully Nonlinear Elliptic Equations} }
 \author{{Nikolai Nadirashvili\thanks{LATP, CMI, 39, rue F. Joliot-Curie, 13453
Marseille  FRANCE, nicolas@cmi.univ-mrs.fr},\hskip .4 cm Serge
Vl\u adu\c t\thanks{IML, Luminy, case 907, 13288 Marseille Cedex
FRANCE, vladut@iml.univ-mrs.fr} }}

\date{}
\maketitle

\def\n{\hfill\break} \def\al{\alpha} \def\be{\beta} \def\ga{\gamma} \def\Ga{\Gamma}
\def\om{\omega} \def\Om{\Omega} \def\ka{\kappa} \def\lm{\lambda} \def\Lm{\Lambda}
\def\dl{\delta} \def\Dl{\Delta} \def\vph{\varphi} \def\vep{\varepsilon} \def\th{\theta}
\def\Th{\Theta} \def\vth{\vartheta} \def\sg{\sigma} \def\Sg{\Sigma}
\def\bendproof{$\hfill \blacksquare$} \def\wendproof{$\hfill \square$}
\def\holim{\mathop{\rm holim}} \def\span{{\rm span}} \def\mod{{\rm mod}}
\def\rank{{\rm rank}} \def\bsl{{\backslash}}
\def\il{\int\limits} \def\pt{{\partial}} \def\lra{{\longrightarrow}}

\section{Introduction}
\bigskip

In this paper we study a class of fully nonlinear second-order
elliptic equations of the form
$$F(D^2u)=0\leqno(1)$$
defined in a domain of ${\bf R}^n$. Here $D^2u$ denotes the
Hessian of the function $u$. We assume that
 $F$ is a Lipschitz  function defined on 
$S^2({\bf R}^n)$ of the space of ${n\times n}$ symmetric
matrices. Recall  that (1) is called
uniformly elliptic if  there exists a constant $C=C(F)\ge 1$
(called an {\it ellipticity constant\/}) such that
$$C^{-1}||N||\le F(M+N)-F(M) \le C||N||\; \leqno(2)$$ for any
non-negative definite symmetric matrix $N$; if $F\in C^1(D)$ then
this condition is equivalent to
$${1\over C'}|\xi|^2\le
F_{u_{ij}}\xi_i\xi_j\le C' |\xi |^2\;,\forall\xi\in {\bf
R}^n\;.\leqno(2')$$

 Here, $u_{ij}$ denotes the partial derivative
$\pt^2 u/\pt x_i\pt x_j$. A function $u$ is called a {\it
classical\/} solution of (1) if $u\in C^2(\Om)$ and $u$ satisfies
(1).  Actually, any classical solution of (1) is a smooth
($C^{\alpha +3}$) solution, provided that $F$ is a smooth
$(C^\alpha )$ function of its arguments.

For a matrix $S \in   S^2({\bf R}^n)$  we denote by $\lambda(S)=\{
\lambda_i : \lambda_1\leq...\leq\lambda_n\}
 \in {\bf R}^n$  the (ordered) set  of eigenvalues of the matrix $S$. Equation
(1) is called a Hessian equation  ([T1],[T2] cf. [CNS]) if the
function $F(S)$ depends only
  on the eigenvalues $\lambda(S)$ of the matrix $S$, i.e., if
 $$F(S)=f(\lambda(S)),$$
 for some function $f$  on ${\bf R}^n$ invariant under  permutations of
 the coordinates.

 In other words the equation (1) is called Hessian if it is invariant under
 the action of the group
 $O(n)$ on $S^2({\bf R}^n)$:
 $$\forall O\in O(n),\; F({^t O}\cdot S\cdot O)=F(S) \;. $$

\medskip
 
Consider the Dirichlet problem
$$\cases{F(D^2u)=0 &in $\Om$\cr
u=\vph &on $\pt\Om\;,$\cr}\leqno(3)$$ where  $\Omega \subset {\bf
R}^n$ is a bounded domain with smooth boundary $\partial \Omega$
and $\vph$ is a continuous function on $\pt\Om$.
   
\medskip

The main goal of this paper is to show that the axially symmetric
solutions of the Dirichlet problem are classical   for Hessian
 elliptic equations. Recall that without the symmetricity assumption this can be
 false in higher dimensions [NV1, NV2].

Let  $\Omega \subset {\bf R}^3$ be a smooth bounded axially symmetric domain.
We consider the Dirichlet problem $(3)$ in $\Omega. $

\medskip
{\bf Theorem 1. }{ \it Let $F\in C^1$  be a uniformly elliptic operator. Let  $\vph \in C^{1,\epsilon }(\partial \Omega )$ be an axially symmetric function, $0<\epsilon <\epsilon_o$, where $\epsilon_o>0$ depends
on the ellipticity constant of $F$.
Then the Drichlet problem $(3)$ has a unique classical
solution $u\in C^2(\Omega )\cap C^{1,\epsilon }(\bar \Omega )$. }
\medskip

{\bf Remark. } The same results hold for the solutions of the $n$-dimensional axially symmetric problems (i.e., for the solutions of the form  $u(x)=u(x_1, x_2^2+...+x_n^2)$ ). 

\medskip
The axially symmetric problems are essentially 2-dimensional.
Outside the axis of symmetry one can rewrite the equations as two-dimensional
fully nonlinear equations with lower order terms. However on the axis of symmetry
the equations became singular and that limits the application of the strong
methods known for the dimension 2.

\section{ Proof of Theorem 1}

Let $\Omega \subset {\bf R}^n$. Let 
$$Lw=  \sum a_{ij}(x) {\partial^2 w \over  \partial x_i \partial x_j }, \leqno (2.1)$$
be a linear uniformly elliptic operator defined in a domain $\Omega \subset {\bf R}^n $,
$$C^{-1}|\xi  |^2 \leq \sum a_{ij}\xi_i\xi_j \leq C|\xi |^2.$$ 

We will need the following  propositions, see [GT], [K].

\medskip
{\bf Proposition 1. } {\it Let $G\subset {\bf R}^n $ be be a bounded domain with a smooth boundary.
Let $u\in C^2 (\bar G) $ be a solution of the equation
$$Lu = 0  \quad  in  \quad   G,$$
$u_{|\partial G} =\phi $. Then
$$||Du||_{C^{\alpha } (\partial G )} \leq C||\vph ||_{C^{1,\alpha }(\partial G )},$$
where positive constants $\alpha $ and $C$ depend on $G$ and the ellipticity constant
of the operator $L$. }

\medskip
{\bf Proposition 2. }{\it Assume that $F\in C^1$, $F(0)=0$, $\partial \Omega \in C^2$ and  the uniform ellipticity condition
$(2')$ holds. 
Let  $u \in C^2(\bar \Omega )$ be a solution of the Drichlet problem $(3).$
Then
$$||Du||_{C^{\alpha } (\partial \Omega )} \leq C||\vph ||_{C^{1,\alpha }(\partial \Omega )},$$
where positive constants $\alpha $ and $C$ depend on  $\Omega $ and on the ellipticity constants
of $F$.}

\medskip
Two following propositions are essentially two-dimensional, see [BJS], [GT].

\medskip
{\bf Proposition 3. } {\it  Let $u\in C^2(D_1)$, where $D_r\subset {\bf R}^2$ be the disk $|x|<r$,
 and let $u$   be a solution in $D_1$ of the equation
 $$Lu=0,$$
 where $L$ is the elliptic operator $(2.1).$ Then
 $$ osc_{D_1}\ u_{x_1} \geq (1+ \xi )osc_{D_{1/2}}\ u_{x_1},$$
 where $\xi>0$ be a constant depending only on the ellipticity constant of operator $L$. }
 
\medskip
 {\bf Proposition 4. }{\it Let  $u\in C^2(D_1)$ be a solution of a fully nonlinear
 elliptic equation
 $$H(D^2u, Du, x)=0$$
 in $D_1$, and $H(0,0,x)=0$. Let $| u|<M$. Then
 $$||u||_{C^{2,\alpha }(D_{1/2})}<CM,$$
 where $\alpha ,C>0$ are constants depending on the ellipticity constant of $H$ and $C^1$-norm of the function $H$. }

\medskip
As a corollary of Proposition 3 we have 

\medskip
{\bf Lemma 1. }{\it  Let $u\in C^2(D_1)$ be a solution of the equation
 $$Lu=0,$$
 in $D_1$ and $l$ an affine  linear function  in $D_1$. Let  $|l-u|<M$. Then for any $\epsilon >0$ there are $\alpha ,r>0$ depending only on $\epsilon$ and the ellipticity constant of $L$ such that  
$$|| u-l||_{C^{1,\alpha }(D_r)} < \epsilon M.$$   }

\medskip
Applying Lemma 1 to the derivative of the solutions of fully nonlinear elliptic equation we get

\medskip
{\bf Lemma 2.} {\it  Let $u\in C^2(D_1)$ be a solution of the fully nonlinear equation
 $$F(D^2u)=0,$$
 in $D_1$ and $F(0)=0$. Let $q$ be a quadratic polynomial  in $D_1$ such that $|q- u|<M$. Then for any $\epsilon >0$ there are $\alpha ,\rho>0$ depending only on $\epsilon$ and the ellipticity constant of $F$ such that  
$$|| u-q||_{C^{2,\alpha }(D_\rho)} < \epsilon M.$$   }
 
\medskip
Proving Theorem 1 we may assume without loss that $F(0)=0$.

\medskip
Let $x_1,x_2,x_3$ be an orthonormal coordinate system in ${\bf R}^3$ and $x_1$ be an
axis of symmetry of the domain $\Omega $. Denote
$$\omega = \{x\in \Omega , x_3=0\}.$$

Let $u$ be a classical axially symmetric solution of the Dirichlet problem (3). Denote
$$||u||_{C(\Omega )} =A.$$

Since $u_{x_3}$ is a solution of  linear uniformly
elliptic equation $Lu_{x_3}=0$ 
 and $u_{x_3}=0$ on $\omega $ then by Proposition 1
$$||u_{x_3x_3}||_{C^{\alpha }(\omega' )} \leq C||\vph ||_{C^{1,\alpha }(\partial G )},\leqno (2.1)$$
where $\omega' \subset \subset \omega $, positive constants $\alpha $ and $C$ depend on $G, \omega'$ and the ellipticity constant
of the operator $F$.

We define two-dmensional Hessian elliptic operators $f_a, a\in R$,
$$f_a(\lambda_1, \lambda_2)= f(\lambda_1, \lambda_2, a ).$$

Let $y\in \Omega $ be a point on the axis $x_1$.
Denote
$$h=dist (y, \partial \Omega ).$$
Define  for  $ 0<r<h$ the function $u_r$ on the unit disk
$D_1\subset R^2$  by
$$u_r(x)=u_r(x_1,x_2)=(u(r(x_1-y_1,x_2))-u(y))/r^2.$$

Set $a=u_{x_2x_2}(y)$.
Let $v_r$ be a solution of the Dirichlet problem
$$\cases{f_a(\lambda (D^2v_r))=0 &in $D_1 $\cr
v_r=u_r &on $\pt D_1\;,$\cr}\leqno(2.2)$$ 
The  classical solution of two-dimensional Dirichlet problem (2.2) is known to
exist, e.g. [GT].

Since  our  equation $F(D^2u)= 0$ is homogeneous 
we can assume without loss that  the  inequalities 
$$1< |\nabla F | < C$$
 hold for a positive constant $C.$

From (2.1) and the last inequalities it  follows   easily that the functions 

\noindent $u_r - C_or^{\alpha }(1-|x|^2) $ and
$u_r + C_or^{\alpha }(1-|x|^2) $ are, for a sufficiently large constant $C_o$, sub- and supersolutions of the Dirichlet problem (2.2).
Hence
$$|u_r-v_r|\leq C_or^{\alpha }.$$

Denote, $w_r=u_r-v_r$. 

Let  $\rho $ be the constant of Lemma 2 for the elliptic operator $f$ and $\epsilon =1/2$. 

Define a sequence of functions $u_n$ in $D_1$, $n=1,2,...$, by
$$u_n=u_{h\rho^n}.$$

Correspondingly we define $v_n=v_{h\rho^n}$, $w_n =u_n-v_n$.

From Lemma 2 we get the following recurrence inequalities:
 there are quadratic polynomials
$q_n$, $n=1,2,...$,  such that $f(q_n)=0$ and 
$$||v_{n+1} -q_{n+1}||_{|C(D_1)} \leq {1\over 2}||v_n -q_n||_{|C(D_1)} + C_o\rho^{\alpha n}.$$
 
Since $|u_1|<A/h^2$, we get
$$||v_n -q_n||_{|C(D_1)}<2AC_o\rho^{\alpha n}/h^2,$$
$$||u_n -q_n||_{|C(D_1)}<2AC_o\rho^{\alpha n}/h^2,$$
$$||w_n||_{|C(D_1)}<C_o\rho^{\alpha n}.$$
for all $n=1,2,...$.

Hence, since the functions $u_n$ are obtained as dilations of $u$, it follows that 
$$||q_{n+1} -q_n||_{|C(D_1)}<2AC_o\rho^{\alpha n-2}/h^2.\leqno (2.3)$$

Therefore 
$$||u_n||<AC_1/h^2\leqno (2.4)$$
for   a constant $C_1>0$ depending only on the ellipticity constant of $F$, $n=1,2,...$.

Denote 
$$E=\{ z=x+y: |x|<h/2, x_2/x_1>1/4 \},$$
$$G= \{x\in D_1 : x_2> 1/4, dist (x, \partial D_1 )> 1/4 \} ,$$
$$G_n=\{ x: x/h\rho^n \in G\},$$
$n=1,2,...$.

Set $g_n=u_n-q_n$. Then from (2.3), (2.4) and Proposition 4 we have
$$||g_n||_{C^{2,\alpha }(G)}<AC_2/h^2,$$
where $C_2>0$ depends only on the ellipticity constant of $F$. Since
$$||g_n||_{|C(D_1)}<2AC_o\rho^{\alpha n}/h^2$$
then by interpolation between the last two inequalities we get
$$||g_n||_{C^{1,\alpha /2 }(G)}<AC_3\rho^{\alpha n/2}/h^2,$$
where $C_3>0$ depends on the ellipticity constant of the equation. Thus 
$$||u||_{C^{1,\alpha /2 }(G_n)}<AC_3/h^2,$$
for all $n=1,2,...$. Together with (2.3) the last inequality gives
$$||u||_{C^{1,\alpha /2 }(E)}<AC_4/h^2,\leqno (2.5)$$
where $C_4>0$ depends only on the ellipticity constant of the equation.

By (2.1) on the axis $x_1$ the second derivatives $u_{x_2x_2}=u_{x_3x_3}$ 
satisfy the H\"{o}lder estimates. Since on the axis the mixed derivatives $u_{x_ix_j}=0$ for $i\neq j$
we conclude from the equation that the second derivative $u_{x_1x_1} $
satisfies the H\"{o}lder estimates as well. These estimates together with (2.5)
give the following inequality
$$ ||u||_{C^{1,\alpha /2 }(\omega')}<AC_5, $$
where $C_5>0$ depends on the ellipticity constant of the equation and the 
distance of $\omega'$ to the boundary $\partial \omega $.

Combining the last inequality with Proposition 2 we get the following apriori  
estimate for the axially symmetric solutions of  fully nonlinear uniformly elliptic equations:
 
\bigskip 
{\bf Lemma 3 }. 
{\it Let $u\in C^2 (\Omega) $ be an axially symmetric solution of $(3)$ and let $\Omega' $ 
 be a compact subdomain of $\Omega $. Then the following inequalities hold:
$$||u||_{C^{1,\alpha } ( \Omega )} \leq C||\vph ||_{C^{1,\alpha }(\partial \Omega )},$$
$$||u||_{C^{2,\alpha } ( \Omega' )} \leq C'||\vph ||_{C^{1,\alpha }(\partial \Omega )},$$
where positive constants $C,C'$ and $\alpha $ depend on  $\Omega $ and on the ellipticity constant of $F$, $C'$ depending also on the distance of $\Omega'$ to the boundary $\partial \Omega $.}

 \bigskip
 The apriori estimate of  Lemma 3  and the standard method of continuation
by parameter, see, e.g., [GT],  gives  the classical solvability of the Dirichlet
problem  (3) for a uniformly elliptic equation.

\bigskip
 \centerline{REFERENCES}

\medskip
 \noindent [CC] L. Caffarelli, X. Cabre, {\it Fully Nonlinear Elliptic
Equations}, Amer. Math. Soc., Providence, R.I., 1995.

\medskip
 \noindent [BJS] L.Bers, F.John, M.Schechter, {\it Partial Differential Equations},
 Interscience Publisher, New York-london-Sydney, 1964.

\medskip
 \noindent [CIL]  M.G. Crandall, H. Ishii, P-L. Lions, {\it User's
guide to viscosity solutions of second order partial differential
equations,} Bull. Amer. Math. Soc. (N.S.), 27(1) (1992), 1--67.

\medskip
 \noindent [CNS] L. Caffarelli, L. Nirenberg, J. Spruck, {\it The Dirichlet
 problem for nonlinear second order elliptic equations III. Functions
  of the eigenvalues of the Hessian, } Acta Math.
   155 (1985), no. 3-4, 261--301.

\medskip
 \noindent [GT] D. Gilbarg, N. Trudinger, {\it Elliptic Partial
Differential Equations of Second Order, 2nd ed.}, Springer-Verlag,
Berlin-Heidelberg-New York-Tokyo, 1983.

\medskip

\noindent [K] N.V. Krylov, {\it Nonlinear Elliptic and Parabolic
Equations of Second Order}, Reidel, 1987.

\medskip
\noindent [NV1] N. Nadirashvili, S. Vl\u adu\c t, {\it On Hessian
fully nonlinear elliptic equations}, arXiv:0805.2694 [math.AP],
submitted.

\medskip
\noindent [NV2] N. Nadirashvili, S. Vl\u adu\c t, 
{\it Nonclassical   Solutions of Fully Nonlinear Elliptic
Equations II: Hessian  Equations  and Octonions  }, arXiv:0912.312 
[math.AP], submitted.

\medskip
\noindent [T1] N. Trudinger, {\it Weak solutions of Hessian
equations,} Comm. Partial Differential Equations 22 (1997), no.
7-8, 1251--1261.

\medskip
\noindent [T2] N. Trudinger, {\it On the Dirichlet problem for
Hessian equations,} Acta Math.
   175 (1995), no. 2, 151--164.

\end{document}